\documentclass[12pt]{amsart}

\usepackage{amsfonts,amsmath,amssymb,amsthm,latexsym,verbatim,amscd}
\usepackage{epsfig, graphicx, psfrag, subfigure}
\usepackage{graphics}
\usepackage{pdfpages}
\usepackage{bbm}
\usepackage{changes}

\newcommand{\calC}{\mathcal{C}}

\newcommand{\calE}{\mathcal{E}}

\theoremstyle{plain}
\numberwithin{equation}{section}

\def \l{\lambda}

\def \o{\omega}

\newcommand{\R}{{\mathbb R}}

\newcommand{\SY}{{\mathbb S }}
\newcommand{\HY}{{\mathbb H}}

\newcommand{\C}{{\mathbb C}}

\theoremstyle{plain}
\newtheorem{theorem}{Theorem}[section]
\newtheorem{theorems}{Slab Theorem}

\newtheorem*{corollary*}{Corollary}

\newtheorem*{proposition*}{Proposition}

\newtheorem{lemmad}{Dragging Lemma}

\newtheorem*{example*}{Example}

\newtheorem*{definition*}{Definition}

\newtheorem*{notation*}{Notation}
\newtheorem{remark}[theorem]{Remark}
\newtheorem*{remark*}{Remark}

\title[Minimal surfaces in a slab of ${\mathbb H}\times {\mathbb R}$]{Properly immersed minimal surfaces in a slab of ${\mathbb H}\times {\mathbb R}$, ${\mathbb H}$ the hyperbolic plane}

\author{P. Collin}
\address{P.Collin, Institut de math\'ematiques de Toulouse,
Universit\'e Paul Sabatier, 118, route de Narbonne, 31062 Toulouse cedex ,France}
\email{collin@math.ups-tlse.fr   }

\author{L. Hauswirth}
\address{L. Hauswirth, Universit\'e Paris-Est, LAMA (UMR 8050), UPEMLV, UPEC, CNRS, F-77454, Marne-la-Vall\'ee, France}
\email{hauswirth@univ-mlv.fr}

\author{H. Rosenberg}
\address{H. Rosenberg, Instituto Nacional de Matematica Pura e Aplicada (IMPA) Estrada Dona Castorina 110, 22460-320, Rio de Janeiro-RJ, Brazil}
\email{ rosen@impa.br}

\begin{document}

\begin{abstract} 
We prove that the ends of a properly immersed simply or one connected minimal surface in ${\mathbb H}\times {\mathbb R}$ contained in a slab of height less than $\pi$ of ${\mathbb H}\times {\mathbb R}$, are multi-graphs. When such a surface is embedded then the ends are graphs.
When embedded and simply connected, it is an entire graph.
\end{abstract}
\thanks{{\it The authors was partially supported by the ANR-11-IS01-0002 grant.} \today}
\maketitle
{
\section{Introduction}
A fundamental problem in surface theory is to understand surfaces of prescribed curvature in homogenous 3-manifolds.  Simply connected properly embedded surfaces are the simplest to consider and after the compact sphere, the plane is next.
There are some unicity results.  A proper minimal embedding of the plane in $\R^3$ is a flat plane or a helicoid \cite{M-R}.  A proper embedding of the plane in hyperbolic 3-space as a constant mean curvature one surface (a Bryant surface) is a horosphere \cite{C-H-R}. Also, a proper minimal embedding of an annulus $\SY^1 \times \R$ in $\R^3$ is a catenoid \cite{collin} and a proper
minimal embedding of an annulus with boundary,-$\SY^1 \times \R^+$-, is asymptotic to an end of a catenoid, plane
or helicoid \cite{M-R, bernstein-brenner}. This is true for Bryant annular ends in $\HY^3$. A proper constant mean curvature
one embedding of an annulus $\SY^1 \times \R$ in $\HY^3$ is a catenoid cousin. Such an annulus
with compact boundary in $\HY^3$ is asymptotic to an end of a catenoid cousin or a horosphere.

In this paper we consider proper minimal embeddings (and immersions) of the plane and the
annulus in $\HY \times \R$.

Contrary to $\R^3$, there are many such surfaces in $\HY \times \R$.  Given any continuous rectifiable curve $\Gamma \subset \partial _{\infty} (\HY \times \R)$, $\Gamma$ a graph over $\partial_\infty (\HY)$, there is an entire minimal graph asymptotic to $\Gamma$ at infinity \cite{NR,NR2}.  Also if $\Gamma$ is an ideal polygon of $\HY$, there are necessary and sufficient conditions on $\Gamma$ which ensure the existence of a minimal graph over the interior of $\Gamma$, taking values plus and minus infinity on alternate sides of $\Gamma$ \cite{C-R}.  This graph is then a simply connected minimal surface in $\HY \times \R$.  We will see there are many minimal embeddings of the plane that are not graphs (aside from the trivial example of a (geodesic of $\HY) \times \R$; a vertical plane.

In this paper we will give a condition which obliges a properly embedded minimal plane to be an entire graph.  More generally, we will give a condition which obliges properly immersed minimal surfaces of finite topology to have multi-graph ends.

We will see that slabs $S$ of height less than $\pi$ play an important role. There are two reasons for this. First, complete
vertical rotational catenoids exist precisely when their height is less than $\pi$. Secondly, for $h > \pi$, there are vertical
rectangles of height $h$ at infinity; i.e., in $\partial _{\infty} (\HY \times \R)$, that bound simply connected minimal surfaces
$M(h)$ (the rectangle is the asymptotic boundary of $M(h)$). These surfaces are invariant by translation along a horizontal
geodesic and we discuss them in detail later; \cite{Hau,SaEarp-Tou,M-R-R}.

Let $\epsilon >0$ and $S = \{ (p,t) \in \HY \times \R ; |t| \leq (\pi-\epsilon)/2 \}$. We will prove that simply or one-connected
properly immersed minimal surfaces in $\HY \times \R$, $\Sigma \subset S$ have multi-graph ends. In fact, an annular properly
immersed minimal surface in $S$, with compact boundary, has a multi-graph subend.  More precisely we prove:

\begin{theorems}
Let $S \subset \HY \times \R$ be a slab of height $\pi - \epsilon$ for some $\epsilon >0$. 
Assume $\Sigma$ is a properly immersed minimal surface  in $\HY \times \R$, $\Sigma \subset S$.
If $\Sigma$ is simply connected and embedded then $\Sigma$ is an entire graph.

More generally;

\begin{enumerate}

 \item If $\Sigma$ is of finite topology and embedded with one end then $\Sigma$ is simply connected and an entire graph.

\item If $\Sigma$ is homeomorphic to $\SY^1 \times \R^+$, then an end of $\Sigma$ is a multi-graph.
If this annular surface is embedded then an end is a graph.

\end{enumerate}
In particular, by (2), if $\Sigma$ is of finite topology then each end of $\Sigma$ is a multi-graph.

\end{theorems}

Results of this nature have been obtained by Colding and Minicozzi in their study
of embedded minimal disks in balls of $\R^3$, whose boundary is on the boundary of the ball;
see proposition III. 1.1 of \cite{C-M-1}. We have been inspired here by their ideas; in particular
using foliations by catenoids to control minimal surfaces.

\begin{remark}
We will give an example of an Enneper-type minimal surface in $S$. This is a properly immersed minimal surface
in $\HY \times \R$, $\Sigma \subset S$  that is simply connected and whose end is a $3$-fold covering graph over the complement of a compact disc in $\HY$.
\end{remark}
\begin{remark}
Also there is an example of a properly embedded simply connected minimal surface in a slab
 of height $\pi$ that is not a graph.  We will describe this surface after the proof of the Slab Theorem.   
 Thus $\pi$ is optimal for the slab theorem.
 \end{remark}

\section{The Dragging Lemma}

\begin{lemmad} Let $g:\Sigma \to N$ be a properly immersed minimal surface in a complete
$3$-manifold $N$. Let $A$ be a compact surface (perhaps with boundary)
and $f: A \times [0,1] \to N$ a $\calC^1$-map such that $f(A \times \{ t \})=A(t)$ is a minimal
immersion for $0 \leq t \leq 1$. If $\partial(A(t)) \cap g(\Sigma) = \emptyset$ for $0\leq t \leq 1$
and $A(0) \cap g(\Sigma) \neq \emptyset $, then there is a ${\calC}^1$ path $\gamma (t)$ in $ \Sigma$, such that $g \circ \gamma(t) \in A(t) \cap g( \Sigma) $ for $0 \leq t \leq 1$. Moreover we can prescribe any initial value $g \circ \gamma (0) \in A(0) \cap g(\Sigma)$.
\end{lemmad}

\begin{remark}
To obtain a $\gamma (t)$ satisfying the Dragging lemma that is continuous (not necessarily $\calC ^1$)
it suffices to read the following proof up to (and including) Claim 1.
\end{remark}

\begin{proof}
When there is no chance of confusion we will identify in the following $\Sigma$ and its image $g(\Sigma)$,
$\gamma \subset \Sigma$ and $g\circ \gamma$ in $g (\Sigma) \subset N$. In particular when we consider embeddings
of $\Sigma$ there is no confusion. 

Let $\Sigma (t)= g(\Sigma) \cap A(t)$ and $\Gamma (t)= f^{-1} (\Sigma (t))$, $0 \leq t \leq 1$ the pre-image in $A \times [0,1]$.

When $g: \Sigma \to N$ is an immersion, we consider $p_0 \in g(\Sigma) \cap A(0)$,
and pre-images $z_0 \in g^{-1} (p_0)$ and $(q_0,0) \in f^{-1} (p_0)$.  We will obtain the arc $\gamma (t) \in \Sigma$ in a neighborhood of $z_0$ by a lift of an arc $\eta (t)$ in a neighborhood of  $(q_0,0)$ in $\Gamma ([0,1])$ i.e. $g\circ \gamma (t)=f \circ \eta(t)$.
We will extend the arc continuously by iterating the construction.

\begin{figure}
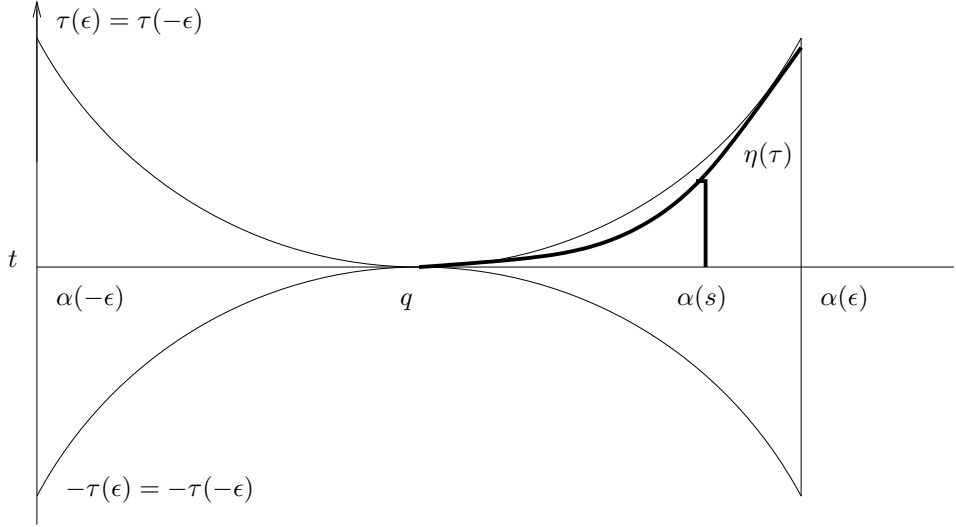
\caption{Neighborhood of a singular point }
\label{fig:figure14d}
\end{figure}

Since $\Gamma (t)$ represents the intersection of two compact minimal surfaces, we know $\Gamma (t)$ is a set of a finite number of  compact analytic curves $\Gamma_1 (t),...,\Gamma_k (t)$. These curves $\Gamma_i (t)$ are analytic immersions of topological circles. By hypothesis, $\Gamma (t) \cap (\partial A \times [0,1])= \emptyset$ for all $t$. The maximum principle assures that  the immersed curves can not contain a small loop, nor an isolated point. Since $A(t)$ is compact and has bounded curvature, a small loop in $\Gamma (t)$ would bound a small disc $D$ in $\Sigma$ with boundary in $A$. Since $A$ is locally a stable surface, we can consider
a local foliation around the disc and find a contradicttion with the maximum principle.
We say in the following that $\Gamma (t)$ does not contain small loops.

{\bf Claim 1:} We will see that for each $t$ with $\Gamma (t) \neq \emptyset$, $t<1$ there is a $\delta (t) >0$ such that if $(q,t) \in \Gamma (t)$, then there is a $\calC ^{1}$ arc $\eta(\tau)$ defined for $t  \leq  \tau \leq t + \delta (t)$ such that $\eta (t) =(q,t)$ and $\eta(\tau) \in \Gamma (\tau)$ for all $\tau$ (there may be values of $t$
where $\gamma ' (t)=0$).

Since $\Gamma (0) \neq \emptyset$, this will show that the set of $t$ for which $\eta (t)$ is defined is a non empty open set. This defines an arc $\gamma (\tau)$ as a lift of  $f \circ \eta (\tau) \subset A(\tau)$ in a neighborhood of $\gamma (t) \in \Sigma$.

First suppose $(q,t)  \in \Gamma (t)$ is a point where  $A(t)=f (A \times \{ t\})$ and $g(\Sigma)$ are transverse
at $f(q,t)$. Let us consider the ${\calC}^1$ immersions
$$F: A \times [0,1] \to N \times [0,1] \hbox{ with } F(q,t)= (f(q,t),t)$$
$$G: \Sigma \times [0,1] \to N \times [0,1] \hbox{ with } G(z,t)= (g(z),t).$$
Let $\hat M=F(A \times [0,1]) \cap G (\Sigma \times [0,1])$ and $M=F^{-1} (\hat M)$ .  $F(A \times [0,1])$ and $G(\Sigma \times [0,1])$ are transverse at 
$p=F(q,t)$. Thus $\hat M$ is a 2-dimensional surface of $N \times [0,1]$ near $p$ .  We consider $X (t)$ a tangent vector field along $\Gamma (t)$ and $JX(t)$ an orthogonal vector field to $X(t)$ in $T_{(q,t)} M$. 
If $\partial / \partial t \perp T_p \hat M$, then $T_p \hat M = T_{f(q,t)} A(t)=T_{f(q,t)} g(\Sigma)$ and 
$(q,t)$ would be a non transverse point of intersection of $A(t)$ and $g(\Sigma)$. 
Thus $<JX(t), \partial / \partial t > \neq 0$ and we can find $\eta (\tau)$  a smooth path, 
defined for $\tau \in [t-\delta (q), t + \delta (q) ]$ such that $\eta (t)=(q,t)$ and $\eta ' (t)=JX(t)$ is transverse to $\Gamma (t)$ at $(q,t)$.

By transversality and $f$ being $\calC ^1$ in the variable $t$, we have  a $\delta (q)>0$ such that for 
$t - \delta (q)  \leq \tau \leq t + \delta (q)$, $A(\tau)$ intersects $f \circ \eta (\tau)$ in a unique point and this point 
varies continuously with $t- \delta (q)  \leq  \tau \leq t + \delta (q)$.
With a fixed initial point in $\Sigma$, a lift of $f \circ \eta (\tau)$, defines
$\gamma (\tau) \in \Sigma$.

Again by transversality, we can find a neighborhood of $(q,t)$ in $\Gamma (t)$
and a $\delta >0 $ so that the above path $\gamma (\tau)$ exists for $t -\delta \leq \tau \leq t + \delta$, through each point in the neighborhood of $q$. It suffices, to look for a local immersion of a neighborhood of $0$ in $T_p M$ into $M$, to obtain a ${\calC}^1$ diffeomorphism $\psi : B(0) \subset T_pM \to M$. $M$ has the structure of a ${\calC}^1$ manifold in a neighborhood of points of transversality and this structure extends to $F^{-1} (M) \subset A\times [0,1]$.

\begin{figure}
\label{fig:figured2}
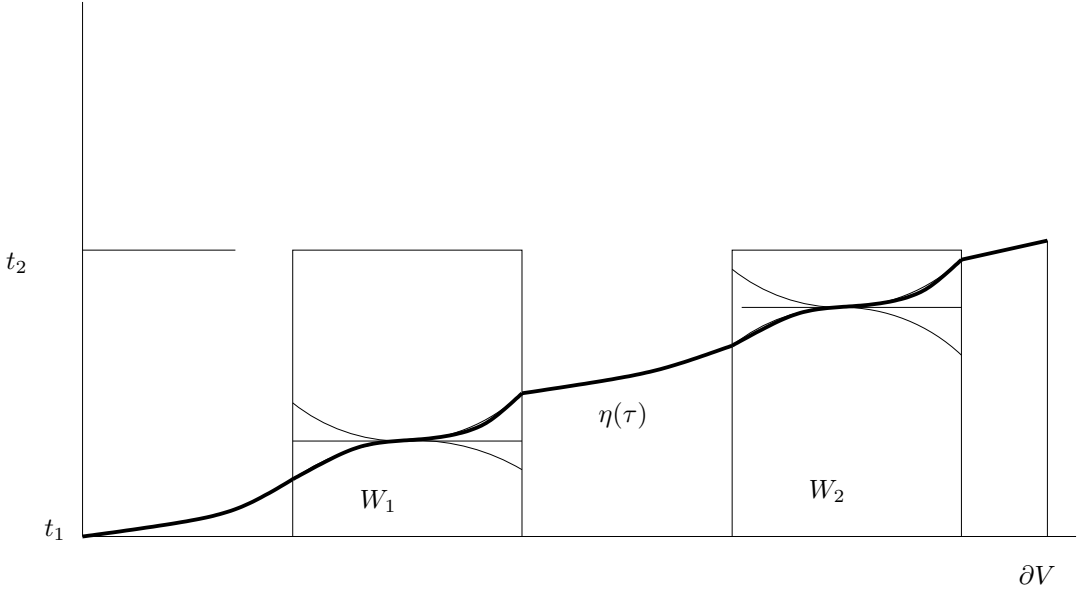
\caption{The curve $\eta (\tau)$ passing through several singularities. }
\end{figure}

We  will find a $\delta >0$ that works in a neighborhood of a singular point $(q,t) \in \Gamma (t)$, where there is a $z \in \Sigma$ such that $f(q,t)=g(z)$ and  $T_{f(q,t)} A(t)=T_{g(z)} g(\Sigma)$. 
We consider singularities of $\Gamma (t)$ where $A(t)$ and $g(\Sigma)$ are tangent. Near a singularity  
$(q,t) \in \Gamma (t)$, $\Gamma (t)$ contains $2k$ analytic curves intersecting at $q$ at equal angles, $k\geq 1$.

Let $V$ be a neighborhood of $q$ in $A$. The set $\Gamma (t) \cap V$ is $2k$ analytic curves. Let $\alpha: ] -\epsilon,\epsilon [ \to V\cap \Gamma (t)$ be a regular parametrization of one curve 
with $\alpha (0) =q$ and $\alpha (\pm \epsilon) \in \partial V$. By transversality as discussed in the previous paragraph $<JX(t), \partial / \partial t > \neq 0$
at $\alpha (s)$ for $s \neq 0$ and $JX (t)$ can be integrated as a curve on $M$ for $t- \delta (s)  \leq \tau \leq t + \delta (s)$.
Here $\delta (s)$ is a ${\calC}^1$ function which can be chosen increasing with $\delta (0)=\delta' (0)=0$.

There exists  a ${\calC}^1$ diffeomorphism  $\phi: \Omega= \{ (s,\tau) \in \R^2; -\epsilon \leq s \leq \epsilon , t - \delta (s)  \leq  \tau \leq t+  \delta (s)\} \to M$ such that $\phi (s,t) = \alpha (s)$ for $s \in ]-\epsilon, \epsilon [$ and $\phi (s ,\tau) \in \Gamma (\tau)$ for  
$t- \delta (s)  \leq \tau \leq t + \delta (s)$. We consider a function $\tau: ]-\epsilon,\epsilon[ \to \R$, such that $(s, \pm \tau (s)) \in \Omega$ and $\tau$ is increasing, $\tau (0)=\tau' (0)=0$ and $\tau (\epsilon)=t+\delta ( \epsilon )$,  $\tau (-\epsilon)=t+\delta (- \epsilon )$.

Now we can construct a path $ \eta (\tau) \in \Gamma (\tau)$ which joins $(q,t)$ to a point in $\Gamma(t+\delta (\epsilon))$. The ${\calC}^1$ arc $f \circ \eta (\tau), t \leq \tau \leq t+\delta (\epsilon)$ is locally parametrized by $\phi (s, \tau (s)), s \in ]0 ,\epsilon[$ and continuously extends to $f(q,t)$ when $\tau \to t$. 
Each point $ \alpha (s)$, can be connected ${\calC}^1$, by the arc
$\phi (s, \tau), t \leq \tau \leq \tau (s)$ from $\alpha (s)$ to  $\phi (s, \tau (s))$, and next a subarc of $\eta (\tau)$ for $\tau (s) \leq \tau \leq t+ \delta (\epsilon)$ (see figure \ref{fig:figure14d}). The constant $\delta (\epsilon)$ depends only on $\alpha (\epsilon)=q_1$, and we note $\delta (q_1):=\delta (\epsilon)$.

Now there are a finite number of arcs $\alpha$ in $V-(q)$, with end points $q$ and a collection
of $q_1,q_2,...q_{2k}$. So one has a $0< \delta $ with $\delta < 
\delta (q_i)$ that works in a neighborhood of $q$. The claim is proved.

To complete the proof of the Dragging Lemma, it suffices to prove that $\gamma (t)$ extends $\calC^1$ for any value of $t \in [0,1]$. Assume that there is a point $t_0$ such that the arc $\gamma (t)$ is defined in a $\calC ^1$ manner for $t < t_0$.
By compactness of $A$, the arc accumulates at a point $(q,t_0) \in \Gamma (t_0)$. Remark that the structure
of $M$ along $\Gamma (t_0)$ gives easily the existence of a continuous extension to $t_0$. To ensure
a ${\calC}^1$ path through $t_0$, we need a more careful analysis at $(q,t_0)$.

{\bf Claim 2:} Suppose the path $\gamma (t)$ satisfies the conditions of the Dragging lemma for
$0 \leq t \leq t_0 <1$. Then $\gamma (t)$ can be extended to $0< t< t_0 + \delta$, to be $\calC^1$ and
satisfy the conditions of the Dragging lemma, for some $\delta >0$.

If $(q,t_0)$ is a transversal point, $M$ has a structure of a manifold and if $t_0-\delta (t_0) < t_1 < t_0$ and $\eta (t_1)=(q_1, t_1)$ is in a neighborhood of $(q,t_0)$, we can find a $\calC ^1$ arc that joins $\eta (t_1)$ to $(q,t_0) \in \Gamma (t_0)$.
Next we extend the arc for $t_0 \leq t \leq t_0 + \delta (t_0)$.

If $(q,t_0)$ is a singular point, we consider a neighborhood $V \subset A$ of $q$ and $\Gamma (t_0)$ intersects $\partial V$ in $2k$
transversal points $q_1,...,q_{2k}$.  We consider $V \times [t_1,t_0]$ with $t_0 - \delta (t_0) <t_1 <t_0$. By transversality at $(q_1,t_0),...,(q_{2k},t_0)$, the analytic set $\Gamma (t_1)$ intersects $\partial V$ in $2k$ points and $V$ in $k$ analytical arcs $\alpha_1,..., \alpha _k$. We suppose that $\eta (t_1) \in \alpha_1 \subset  V \times \{t_1\}$. We construct below a monotonous $\calC ^1$ arc from $\eta (t_1)$ to a point $(\hat q,t_2)$ on $\partial V \times \{t_2\}$ for some $t_1 < t_2 < t_0$ and by transversality
an arc from $(\hat q,t_2)$ to a point $(q',t_0) \in \partial V \times \{t_0\}$, using the fact that $ t_0 - \delta (t_0) < t_2$. Next we can extend the arc in a $\calC^1$ manner from $(q',t_0)$ to some point in $\Gamma (t_0 + \delta (t_0))$.

\begin{figure}
\label{fig:figured3d}
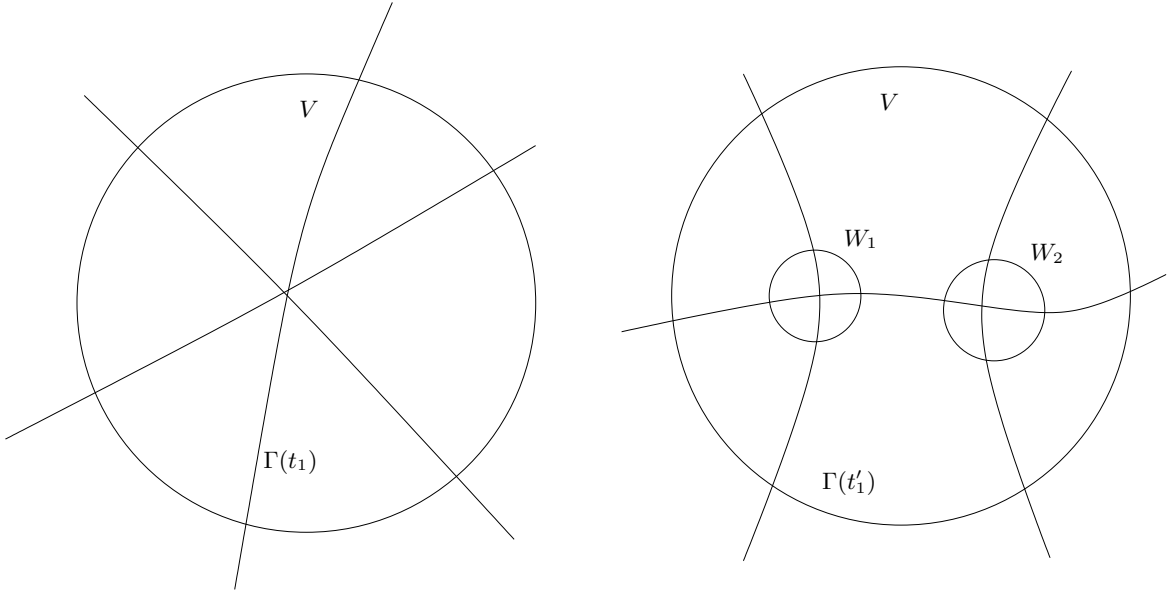
\caption{Left: The curve $\Gamma (t_1)$-Right: The curve $\Gamma (t_1')$.}
\end{figure}

We consider $(\tilde q_1, t_1),...,(\tilde q_\ell, t_1)$ singular points of $\Gamma (t_1) \cap V\times \{t_1\}$ and we denote by $W_1,...,W_\ell$ neighborhoods of $\tilde q_1,...,\tilde q_\ell$ in $A \cap V$. The arc $\alpha_1$ cannot have double points in $V$ without creating small loops. Hence $\alpha_1$ passes through  each $W_1,...,W_\ell$ at most one time, before  joining a point of $\partial V$ (We can restrict $V$ in such a way that there are no 
small loops in $V$). 

First we assume that there is $t_2$ such that for any $t \in [t_1,t_2]$, the curve $\Gamma (t)$ has exactly one isolated singularity in each neighborhood $W_i \times \{t \}$ with the same type
as $\tilde q_i \in \Gamma (t_1)$ ($i=1,...,\ell$) and $t_2< t_1 + \delta (t_1)$. 
If we parametrize $\alpha_1 : [s_0,s_{2\ell+1}]Ê\to \Gamma (t_1)$, we can find $s_1,...,s_{2\ell}$ such that $\alpha_1 (s_{2k-1}), \alpha_1 (s_{2k}) \in  \partial W_k$ and $I_k =[s_{2k-2}, s_{2k-1}]$ are intervals parametrizing transversal points in $\Gamma (t_1)$.

The manifold structure of $M$ gives an immersion $\psi _j: I_j \times [t_1 , t_1 + \delta ] \to M$, $t_1+ \delta < t_2$ and $j=1,...,\ell+1$. In the construction of $\eta$ up to $t_1$, the singular points are isolated;
 then we can assume $\eta (t_1)$ is a regular point of $\Gamma (t_1)$, hence is contained in an
 $\alpha_1 (I_j)$. We construct the beginning of the arc $\eta (\tau)$ as the graph parametrized by
 $\phi_j (s , \tau (s))$ with $\tau$ an increasing function from $t_1$ to $t_1 + \delta/n$ as $s$ varies
 from $\hat s \in I_j$, corresponding to the initial point $\eta (t_1)= \alpha_1 (\hat s)$, to $s_{2j-1}$. Next
 we pass through the singularity $(\tilde q_j, t_1 + 2 \delta/n)$ by constructing
an arc wich joins the point $\phi_j (s_{2j-1},t_1+\delta/n) \in \Gamma (t_1 + \delta/n) \cap \partial W_j$ to
the point $\phi _{j+1}(s_{2j}, t_1 + 3\delta/n) \in \Gamma (t_1 + 3\delta/n) \cap \partial W_j$ (see figure \ref{fig:figured2}).
For a suitable value of $n$ we can iterate this construction, passing through the singularities
$\tilde q_j, \tilde q_{j+1}...$, until we join a point $(\hat q, t_2)$ of $\partial V \times \{ t_2 \}$ and then
we extend the arc up to $t_0$ by transversality outside $V$.

Now we look for this interval $[t_1,t_2]$. Let $t_1<t'_1<t_0$ and $\Gamma (t'_1)$ have several singularities in
some neighborhood $W_k$, or a unique singularity of index less the one of the $\tilde q_k$. We consider in this $W_k$ a finite collection of neighborhoods of isolated singularities $W'_{k,1},...W'_{k,\ell'}$.  We observe, by transversality that there are the same number of components of $\Gamma (t_1)$ and $\Gamma (t'_1)$ in $W_k$  (see figure 3). Hence each $W'_{k,j}$ contains a number of components of  $\Gamma (t'_1)$ strictly less than the number of components of $\Gamma (t_1)$ in $W_k$. The index of the singularity is strictly decreasing along this procedure. We can iterate this analysis up to a point
where each singularity can not be reduced to a simple one. This gives the interval $[t_1,t_2]$.
\end{proof}

\section{Proof of The Slab theorem}

Assume $\Sigma$ is either simply connected ($\partial \Sigma = \emptyset$) or $\Sigma$ is an annulus with compact boundary: $\Sigma$ homeomorphic to $\SY^1 \times \R^+$. Also assume $\Sigma$ is properly minimally immersed in $\HY \times \R$ and $\Sigma \subset S = \{ (p,t) \in \HY \times \R ; | t|< (\pi-\epsilon)/2 \ \}$. 

We fix $h \geq 2 \pi$ sufficiently large so that there are points of $\Sigma$ in the geodesic ball $B$ of radius $h$ of $\HY \times \R$ with center at a point $p_0$ in $\HY \times \{0\}$, and  $\partial \Sigma \subset B$ (if $\partial \Sigma \neq \emptyset$). 

Let ${\rm Cat}(p_0) $ denote a compact part of a rotational catenoid, a bi-graph over $t = 0$, bounded by two circles outside of the slab $S$, $p_0 \in {\rm Cat}(p_0) $.  We assume $h$ is suficientlly large so that $B$ also contains ${\rm Cat}(p_0) $; see figure \ref{fig:figure2bis}.

Since $\Sigma$ is properly immersed, the set $B \cap \Sigma$ has a finite number of compact connected components and there is a compact
 $K$  in $\HY \times \R$, such that any two points of $\Sigma$ in $B$ can be joined by a path of $\Sigma$ in $K$. 

Now suppose $p \in \Sigma - K$ has a vertical tangent plane and let $M= \alpha \times \R$ be tangent to $\Sigma$ at $p$, $\alpha$ a geodesic of $\HY$.
We will prove $M$ must intersect $K$. Suppose this were not the case. $M$ separates $\HY \times \R$ in two components
$M(+)$ and $M(-)$; assume $K \subset M(+)$. The local intersection of $M$ and $\Sigma$ near $p$ consists of $2k$ curves through
$p$, $k \geq 1$, meeting at equal angles $\pi/k$ at $p$. 

\begin{figure}
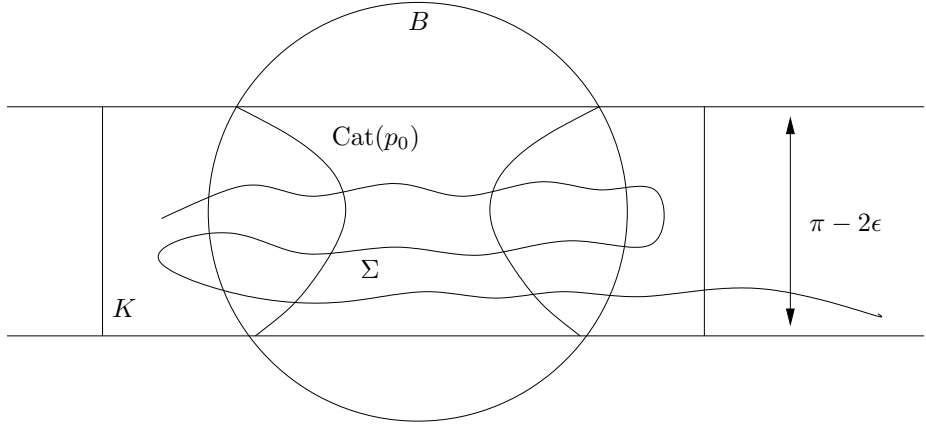
\caption{The compact $K$ and the catenoid ${\rm Cat}(p_0)$}
\label{fig:figure2bis}
\end{figure}

Let $\Sigma_1 (+)$ and $\Sigma _2 (+)$ be distinct local components at $p$ of $\Sigma -  M$, that are contained in
$M(+)$. Observe that $\Sigma_1 (+)$ and $\Sigma_2 (+)$ are contained in distinct components of $\Sigma \cap M(+)$. Otherwise we could find a path $\alpha_0$ in $\Sigma \cap M(+)$, joining a point $x \in \Sigma_1 (+)$ to a point $y \in  \Sigma_2 (+)$. Then join  $x$ to $y$ by a local path $\beta_0$ in $\Sigma$ going through $p$, but $\beta_0 \subset M(+)$ except at $p$; see figure \ref{fig:figure2}.

Let $\Gamma= \alpha_0 \cup \beta_0 \subset \overline{M(+)}$. If $\Sigma$ is simply connected $\Gamma$ bounds a compact disk $D$
in $\Sigma$. If $\Sigma$ is an annulus there are two cases.  $\Gamma \cup \partial \Sigma$ bounds an immersed compact
annulus (we also call $D$) or $\Gamma$ bounds a compact disk in the annulus . By construction of $\Gamma$, $D$ contains points in $M(-)$. 
But $D$ is compact and minimal, $\partial D \subset \overline{M(+)}$, so there would be an interior point of $D$
that is furthest away from $M$ in $M(-)$. This contradicts the maximum principle.

Thus $\Sigma_1 (+)$ and $\Sigma_2 (+)$ are in distinct components $\Sigma_1$ and $\Sigma _2$  of $\Sigma \cap M(+)$. Now let $\mu (\epsilon)$ be the geodesic
of length $\epsilon$ starting at $p$, normal to $M$ at $p$, and contained in $M(+)$.  We will now also denote by $\alpha$, the geodesic $\alpha$  translated vertically to pass through $p$.   Let $\alpha (\epsilon)$ be the complete
geodesic of $\HY$ obtained from $\alpha$ by translating $\alpha$ along $\mu (\epsilon)$ to the endpoint of $\mu (\epsilon)$ distinct
from $p$. The distance between $\alpha$ and $\alpha (\epsilon)$ diverges at infinity; the two geodesics have distinct
end points at infinity.

Choose $\epsilon$ small so that $M(\epsilon)=\alpha (\epsilon) \times \R$ meets both $\Sigma_1 (+)$ and
$\Sigma_2 (+)$, at points $x \in \Sigma_1, y \in \Sigma_2$. If $\Sigma$ is simply connected then no connected component
of $M(\epsilon) \cap \Sigma$ can be compact. Hence $x$ and $y$ are in non-compact components
$C(x) \subset \Sigma_1 \cap M(\epsilon)$ and $C(y) \subset \Sigma_2 \cap M(\epsilon)$.

We claim this also holds when $\Sigma$ is an annulus; more precisely:

{\bf Claim:} If $\Sigma$ is an annulus then the components $C(x)$ of $x$ in $\Sigma_1 \cap M(\epsilon)$ and
$C(y)$ of $y$ in $\Sigma_2 \cap M(\epsilon)$ are both non compact.

{\it Proof of the Claim.} First suppose $C(x)$ and $C(y)$ are compact. Neither can be null homotopic
in $\Sigma$ by the maximum principle, hence $C(x) \cup C(y)$ bound an immersed compact annulus $D$ in 
$\Sigma$. But $\partial D \subset M(\epsilon)$ which is impossible.

It remains to show one of the $C(x), C(y)$ can not be non compact. Suppose $C(x)$ is not compact and $C(y)$ is compact.
So $C(y) \cup \partial \Sigma$ bound an immersed annulus $D$ in $M(\epsilon)(+)$. The distance between
$M$ and $M(\epsilon)$ diverges and $\Sigma$ is proper so there are points $z$ of $C(x)$ arbitrarily far from
$M$. 

Choose such a $z$ so that one can place  a vertical catenoid ${\rm Cat}(z)$ in $M(+)$ which
is a horizontal translation of  ${\rm Cat}(p_0)$ and contains $z$.
The boundary $\partial({\rm Cat} (z)) \cap S = \emptyset$. Let $\eta$ be a geodesic joining
$z$ to the point $p_0 \in B$. Apply the Dragging Lemma to the translation of ${\rm Cat}(z)$ along $\eta$ from $z$ to $p_0$, so that  the translation of ${\rm Cat}(z) $  is contained in $B$, at the end of the movement.
 \begin{figure}
\includegraphics[height=2.5in,width=3.5in]{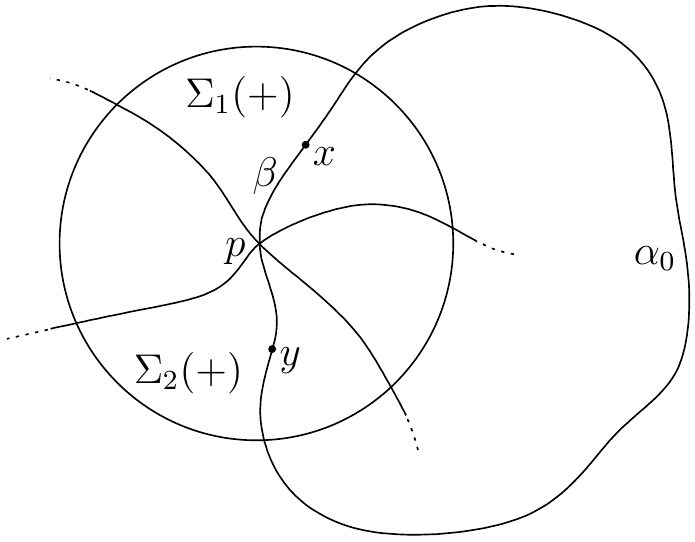}
\caption{A loop $\Gamma$ }
\label{fig:figure2}
\end{figure}

We obtain a path in $\Sigma \cap M(+)$ joining $z$ to a point $\o$ of $\Sigma$ in $B$. Then join $\o$ to a point of $\partial \Sigma$
in $\Sigma \cap K$. This path contradicts that $\Sigma_1$ and $\Sigma_2$ are in distinct components of $\Sigma -M$. Hence the claim
is proved and both $C(x) \subset \Sigma_1 \cap M(\epsilon)$ and $C(y) \subset \Sigma_2 \cap M(\epsilon)$ are non compact.

We now continue the proof that $M$ must intersect $K$. $\Sigma$ is either simply connected or an annulus
and we have two non compact components $C(x) \subset \Sigma_1 \cap M(\epsilon), C(y) \subset \Sigma_2 \cap M(\epsilon)$ and $\Sigma_1,\Sigma_2$ are distinct components of $\Sigma -M$.

Now, as we argued in the proof of the claim above , we can take points $z_1 \in C(x), z_2 \in C(y)$, far enough away from $M$ so that one can place compact vertical
catenoids ${\rm Cat}(z_1)$ and ${\rm Cat}(z_2)$, so $z_1 \in {\rm Cat}(z_1)$, $z_2 \in {\rm Cat}(z_2)$, $\partial {\rm Cat}(z_i) \cap S = \emptyset, i=1,2$ and ${\rm Cat}(z_i)$ are symmetric with respect to $t=0$. Also ${\rm Cat}(z_i) \subset M(+), i=1,2$.

Take horizontal geodesics $\eta_1, \eta_2 \subset M(\epsilon)(+)$ from ${\rm Cat}(z_1)$ to $p_0$ and from ${\rm Cat}(z_2)$ to $p_0$.  Apply the Dragging Lemma along $\eta_1,\eta_2$, to find a path in $\Sigma$ from $z_1$ to a point $\omega_1 \in \Sigma \cap B$,
and another path from $z_2$ to $\o_2 \in \Sigma \cap B$. Join $\o_1$ to $\o_2$ by a path in $\Sigma \cap K$. This contradicts
that $z_1$ and $z_2$ are in distinct components of $\Sigma - M$. Hence $M$ intersects $K$.

To complete the proof of the theorem we will show that when $p \in \Sigma-K$ is far enough from $K$ then $p$ can not
have a vertical tangent plane. We will do this by showing such a vertical tangent plane can not intersect $K$.

To do this we introduce comparison surfaces $M(h)$, first introduced by Hauswirth \cite{Hau}, then by Toubiana and
Sa Earp and Toubiana \cite{SaEarp-Tou}, Daniel \cite{daniel} and Mazet, Rodriguez, Rosenberg \cite{M-R-R}.
The surfaces $M(h)$ are all congruent in $\HY \times \R$, they are complete minimal surfaces invariant
by hyperbolic translation along a geodesic of $\HY$. They exist for each $h > \pi$; we state the properties
we will use.

\begin{figure}
\resizebox{6cm}{!}{\includegraphics{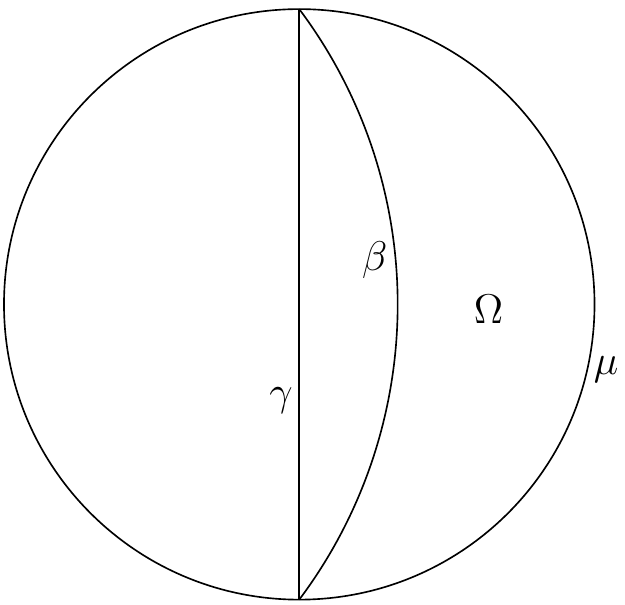}}
\caption{Domain $\Omega$ }
\label{fig:figure3}
\end{figure}

\begin{enumerate}
\item Let $\beta$ be an equidistant of a geodesic $\gamma$ of $\HY \times \{ 0\}$, whose distance
$d_0$ to $\gamma$ is determined by $h$. There is an $M(h) (=M(h,\beta))$ which is a minimal bigraph over the
domain $\Omega$ of $\HY$ indicated in figure \ref{fig:figure3}.
\item $M(h)$ has height $h/2$ over $\Omega \times \{ 0\}$,

\item $M(h)$ meets $\HY \times  \{ 0\}$ orthogonally along $\beta$,

\item $M(h)$ meets each $\HY \times  \{ t\}$ in an equidistant of $\gamma \times \{ t\}$ for $|t| < h/2$,

\item The asymptotic boundary of $M(h)$ is the vertical rectangle of $\HY \times \R$:
$$\partial _{\infty} (M(h))=(\partial _{\infty} (\beta) \times [-h,h]) \cup (\mu \times \{\pm h/2\}).$$
$\mu$ is the arc in $\partial _{\infty} (\HY)$ joining the end points of $\beta$, indicated in figure \ref{fig:figure3}.
\item Given $q \in \HY$ and a vector $v$ tangent to $\HY$ at $q$, there exists an $M(h,\beta)=M(h)$
such that $q \in \beta$ and the tangent to $\beta$ at $q$ is orthogonal to $v$.

\item Now assume $h \geq 2 \pi$ and $S$ is the slab of height $\pi-\epsilon : S =\{ |t| \leq  (\pi-\epsilon)/2\}$.
Then there exists $d_2 >0$ such that if $|t| < (\pi-\epsilon)/2$ and $T_t$ is a vertical translation by $t$
then ${\rm dist}(T_t(M(h)) \cap S, \gamma \times \R) \leq d_2$. Here $\gamma$ is the geodesic
of which $\beta$ is the equidistant. In fact, in any strict subslab $\tilde S$ of $\HY \times ]-h/2,h/2[$, $\tilde S$ symmetric about $t=0$,
$M(h)$ is a bigraph over the moon between $\beta$ and another equidistant $\tilde \beta$ of $\gamma$
and this moon is a bounded distance $d_2$ from $\gamma$.

\item If $\gamma_1$ and $\gamma_2$ are complete geodesics of $\HY$ then $M_{\gamma_1} (h)$ is congruent
to $M_{\gamma_2} (h)$ by a height preserving isometry of $\HY \times\R$.

\item Let $p_0=(0,y_0)$ and for $0<y<y_0$, let $\beta$ be an equidistant curve of $\gamma$ such that $\beta$ is  tangent at $(0,y)$
to the geodesic $\alpha$ through $(0,y)$ and $p_0$
(here $\beta$ is an equidistant of an $M(h)$). Then ${\rm dist}_{\HY} (p_0,\gamma) \to \infty$,
as $y \to 0$; see figure \ref{fig:figure4}.
\end{enumerate}

\begin{figure}
\includegraphics[height=2in,width=2.5in]{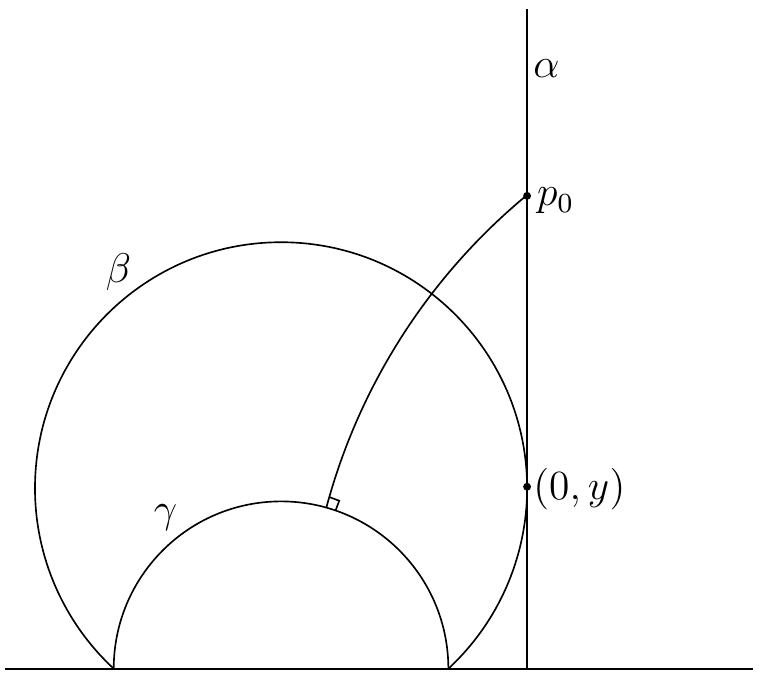}\hfill
\includegraphics[height=2in,width=2.5in]{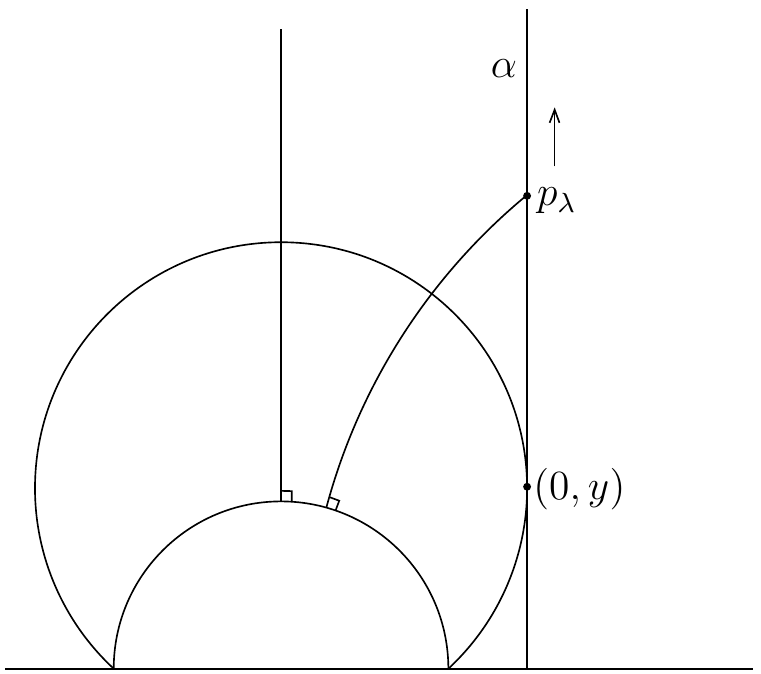}
\caption{Minimizing geodesic from $p$ to $\gamma$}
\label{fig:figure4}
\end{figure}
Here is one way to verify this last assertion. For $(0,y)$ and $h \geq 2 \pi$ fixed, $\beta,\gamma$ and
$M_{\gamma}$ are uniquely determined so that $M_{\gamma}$ is tangent to $\alpha \times \R$ at $(0,y)$ at a point of $\beta$.

So it suffices to fix $(0,y)$, and let $p_{\l} = \l p_0$, with $\l \to \infty$. The minimizing
geodesic from $p_{\l}$ to $\gamma$ tends to the vertical geodesic going up from the point of
maximum $y$ coordinate on $\gamma$; see figure \ref{fig:figure4}. Hence it's length goes to $\infty$ as $\l \to \infty$; see figure \ref{fig:figure4}.

Now we are ready to finish the proof of the theorem. Let $C>0$ be such that if ${\rm dist}_{\HY} (p_0,(0,y)) \geq C$, then ${\rm dist}_{\HY} (p_0,\gamma) > d_2 + {\rm diam} (K)$. Recall that $p_0$ is the center of the ball $B$.
We write $p_0=(0,y_0)$. We will prove that a point of $\Sigma$ at a distance at least $C$ from $p_0$, can not
have a vertical tangent plane; this will prove the theorem.

Assume the contrary. Write $p=(0,y,t)$, $p\in \Sigma$ is a point with a vertical tangent plane and
${\rm dist}_{\HY} (p_0, (0,y)) \geq C$. We vertically translate $M(h)$ by some $t$, $|t|< \pi/2$, so that
$M_t(h) = T_t (M(h))$ is tangent to $\Sigma$ at $p$. Let $\gamma$ be the geodesic associated to
$M_t (h)$. Since $ {\rm dist}_{\HY} (p_0, (0,y)) \geq C$, we know that ${\rm dist}_{\HY} (p_0, \gamma ) > d_2 + {\rm diam} (K)$.

Now if $M_t (h)$ does not intersect $K$, the same proof we gave using $M_t(h)$ in place of $M =\alpha \times \R$, shows we obtain a contradiction.

We explain this further.  Let  $\eta$ be a geodesic of $\HY$ orthogonal to $\gamma$, $p \in \eta$.  The set of all geodesics $\gamma(s)$ , $s \in \R$,  of $\HY$, orthogonal to $\eta$, foliates $\HY$ and $\gamma$ is a leaf of this foliation  For $h > \pi$, the $M(h,\gamma(s))$ foliate the slab of $\HY \times \R$ between $t = h$ and $t = -h$.  Hence if $D$ is a compact minimal surface whose boundary is in some $M(h, \gamma(s))$, then $D$ is contained in $M(h, \gamma(s))$.
The other property of this foliation one uses is the following:  If $\Sigma$ is a properly immersed minimimal surface in the slab $S$ (height $\pi - \epsilon$) and the intersection of $\Sigma$ and $M(h,\gamma(\epsilon))$ has a non compact component $C$, for some $\epsilon > 0$, then ${\rm dist}(C,M(h,\gamma(0)))$ tends to infinity as one diverges in $C$.

Hence it suffices to show $M_t(h)$ can not intersect $K$. Suppose there is a point $w \in M_t (h) \cap K$.
Then ${\rm dist}(\omega, p_0) \leq {\rm diam} K$, and ${\rm dist}(\omega, \gamma) \leq d_2$, so
${\rm dist} (p_0, \gamma) \leq d_2 + {\rm diam }K$; a contradiction.

Now suppose $\Sigma$ of the theorem is embedded. Let $r>0$ and $C(r)=\{ p \in \HY; d(p,p_0)=r \}$. Define
${\rm Cyl} (r)= C(r) \times \R$; ${\rm Cyl} (r)$ is a vertical cylinder of radius $r$.

Let $r_0 >0$ be large so that $\Sigma$ is a multigraph for $r \geq r_0$ (we proved $\Sigma$ is not vertical for large $r$). $\Sigma$ is proper, so $\Sigma \cap {\rm Cyl} (r)$
is a finite union of embedded Jordan curves, each a graph over $C(r)$, for each $r \geq r_0$.

Let $\beta (r)$ be one of the graphical components of $\Sigma \cap {\rm Cyl} (r)$, $r \geq r_0$. If $\Sigma$ is
simply connected (so $\partial \Sigma = \emptyset$),
then the usual proof of Rado's theorem shows $\beta (r)$ bounds a unique compact minimal surface that is a graph over the disk
bounded by $C(r)$. Hence $\Sigma$ is an entire graph.

The same arguments shows that  if $\partial \Sigma = \emptyset$ and $\Sigma$ has an annular end than $\Sigma \cap {\rm Cyl} (r)$
is one Jordan curve, a graph over $C(r)$, that (by Rado's theorem) bounds a unique compact minimal surface, a graph
over the disk bounded by $C(r)$. This proves $(1)$ of the Slab theorem.
\begin{figure}
\resizebox{7cm}{!}{\includegraphics{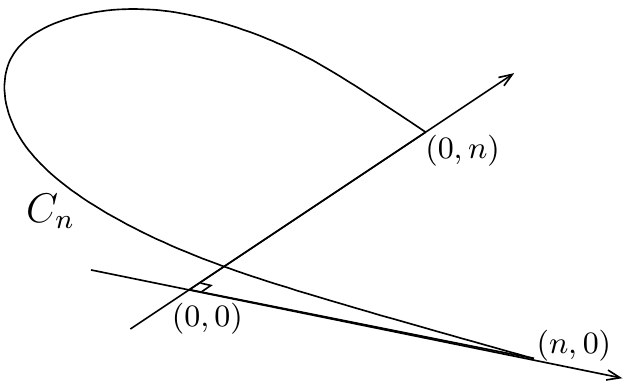}}\hfill
\resizebox{7cm}{!}{\includegraphics{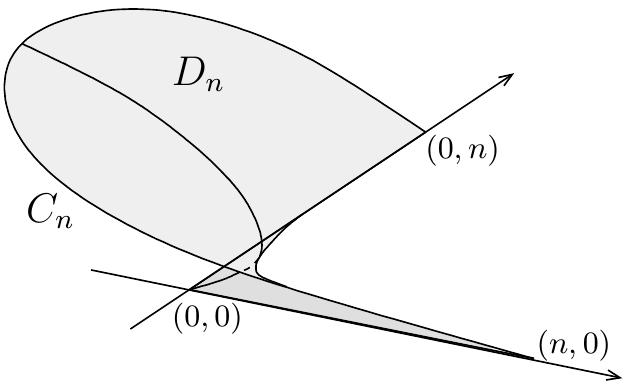}}
\caption{Left: A Jordan curve $C_n$. Right: An Enneper's sector $D_n$}
\label{fig:figure6}
\end{figure}

Now suppose $\Sigma$ is a properly immersed minimal annulus in $S$, with compact boundary. Choose $r$
large so that $\partial \Sigma \subset {\rm Cyl} (r)$ and $\Sigma \cap {\rm Cyl} (r)$ is a finite union of multi-graphs.
Since $\Sigma$ has one end, there is exactly one multi-graph, and the end is a multi-graph.
 If $\Sigma$ is embedded $\Sigma \cap {\rm Cyl} (r)$ is a graph and the end is a graph. This proves (2).

\begin{remark}
We now describe the surface discussed in Remark 1.2.  Consider the surface $M_{\gamma}(h), h > \pi$ and $\gamma$ a geodesic of $\HY$.  Fix a point $p$ in $\gamma$ and deform $\gamma$ through equidistant curves $\beta(t)$, $0 \leq t \leq 1$, through $p$, such that $\beta(0) = \gamma$ and $\beta(1)$ is a horocycle (do this by continuously deforming the endpoints of $\gamma$ to make them converge to one point at infinity).  
Then the surfaces $M_{\beta(t))}(h(t))$ converge to a minimal surface $\Sigma$ of height $\pi$.  The asymptotic boundary of $\Sigma$ is a vertical segment of height $\pi$. $\Sigma$ is vertical along $\beta(1)$, this simply connected embedded minimal surface shows the Slab Theorem fails for slabs of height $\pi$. Benoit Daniel explicitly parametrized this surface; cf proposition 4.17 of \cite{daniel2}.
 \end{remark}

\section{An Enneper type minimal surface}
We construct a properly immersed simply connected minimal surface in $\HY \times \R$ contained in a slab $S$ (of any height $h$)
that is a 3-sheeted multi-graph outside of a compact set.

The idea comes from Enneper's minimal surface $\calE$ in $\R^3$, whose Weierstrass data is $(g,\omega)=(z,dz)$ on the complex plane
$\C$. We think of $\calE$ as constructed by solving Plateau problems for certain Jordan curves in $\R^3$; passing to the limit,
and then reflecting about the two rays in the boundary. More precisely let $C_n$ be a Jordan curve consisting of the segments on
the $x$ and $y$ axes between $0$ and $n$, and the (slightly tilted up) large arc on the circle of radius $n$ (centered at $(0,0)$)
joining $(n,0)$ to $(0,n))$; see figure \ref{fig:figure6}.

Let $D_n$ be a disk of minimal area with $\partial D_n=C_n$; see figure \ref{fig:figure6}-Right.

One can choose $C_n$ so the $D_n$ converge to a minimal surface (not flat) with boundary the positive $x$ and $y$ axes.
$\calE$ is then obtained by Schwarz reflection in the $x$ and $y$ axes. We now do this in $\HY \times \R$.

We use now the unit disk $\{ x^2 + y^2<1\}$ with the hyperbolic metric as a model for $\HY$. Let $h >0$ and $n$ an integer, $n \geq 1$. In $\HY \times \R$,
let $C_n$ be the Jordan curve; see figure \ref{fig:figure8}:

\begin{figure}
\includegraphics[height=2.5in,width=4in]{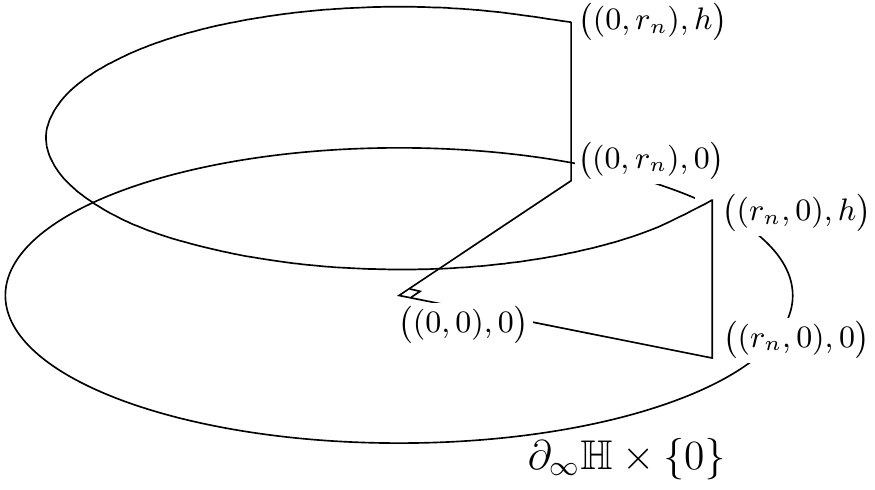}
\caption{A Jordan curve in $\HY \times \R$ }
\label{fig:figure8}
\end{figure}

\begin{figure}
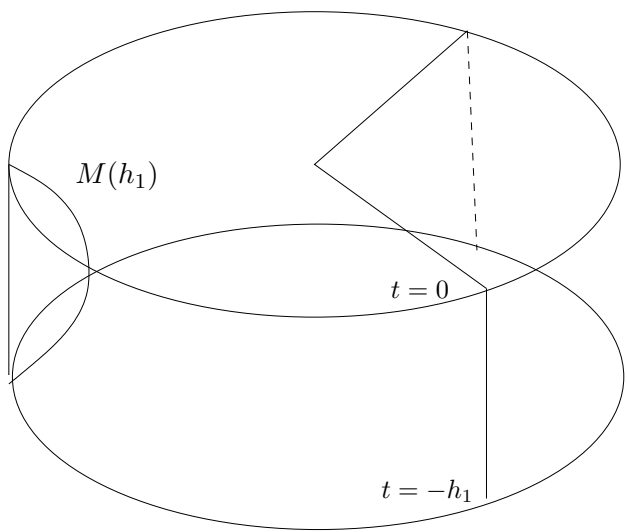
\caption{The barrier $M(h_1)$ }
\label{fig:figure9}
\end{figure}
\[\begin{aligned}
C_n &  = \{ (x,0) ; 0 \leq x \leq r_n \} \cup \{ (0,y) ; 0 \leq y \leq r_n \} \\
&\;\; \cup \{ (r_n \cos \theta, r_n \sin \theta, h) ; \pi/2 \leq \theta \leq 2 \pi \} \\
&\;\; \cup \left( (r_n,0) \times [0,h] \right) \cup \left( (0,r_n) \times [0,h] \right)
\end{aligned}
\]
where $r_n \to 1$.
Let $D_n$ be a least area minimal disk with $\partial D_n = C_n$. We claim a subsequence of $D_n$ converges
to a stable non trivial minimal surface with boundary the positive $x$ and $y$ axes. $\Sigma$ is then obtained by 
reflection in the horizontal boundary geodesics.

First observe the $C_n$ converge in the model of the disk to 

\[\begin{aligned}
C &  = \{ (x,0) ; 0 \leq x \leq 1 \} \cup \{ (0,y) ; 0 \leq y \leq 1 \} \\
&\;\; \cup \{ ( \cos \theta,  \sin \theta, h) ; \pi/2 \leq \theta \leq 2 \pi \} \\
&\;\; \cup \left( (1,0) \times [0,h] \right) \cup \left( (0,1) \times [0,h] \right)
\end{aligned}
\]
The latter circular arc at infinity- $\{ (\cos \theta, \sin \theta,h) ; \pi/2 \leq \theta \leq 2 \pi \}$ is the upper
part of a vertical rectangle at infinity that bounds a minimal surface $M(h_1)$, $h_1>\pi$. We remark that
as $h_1 \to \pi$, $h_1> \pi$ the $\beta$ of $M(h_1)$ tends to a horocycle. Hence $d_2 = {\rm dist} (\beta, \gamma) \to \infty$ as $h_1 \to \pi$. This implies $M(h_1)$ is a barrier for each of the
$D_n$, for $h_1$ sufficiently close to $\pi$; see figure \ref{fig:figure9}.

\begin{figure}
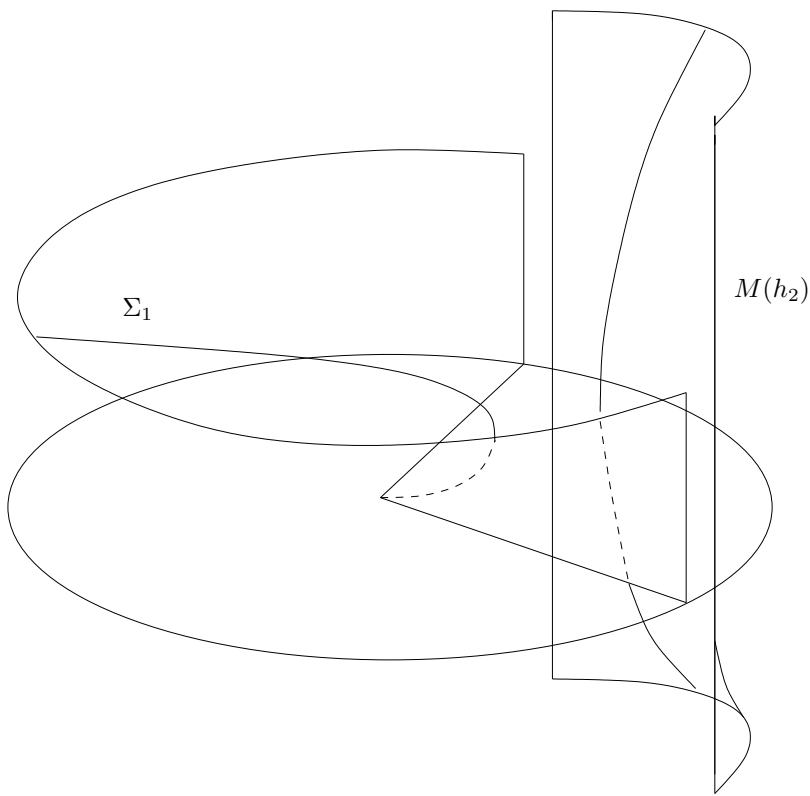
\caption{The barrier $M(h_2)$}
\label{fig:figure10}
\end{figure}

Translate $M(h_1)$ down to be below height $0$, then
go back up to height $h$. There can be no contact with $D_n$. To prevent $D_n$ from escaping in the sector
$0 \leq \theta \leq \pi/2$, one can place another $M(h_1)$ in this sector, whose vertical rectangle at infinity has horizontal arc of
angle less than $\pi/2$ and whose vertical sides go from a height below zero to a height above $h$.

Let $F$ be the vertical solid cylinder in $\HY \times \R$ over the circle of radius $1$ in $\HY$, centered at $(-2,-2)$
(hyperbolic distances). For $n>2$, each $D_n$ intersects $F$ in a minimal surface of bounded curvature
(bound independent of $n$). Hence a subsequence of these surfaces converges to a minimal surface, which 
is not part of a horizontal slice. The barriers $M(h_1)$ then show a subsequence of the $D_n$ converge to a minimal
surface $\Sigma_1$ with boundary the positive $x$ and $y$ axes. 
To guaranty that $D_n$ does not escape in the sector $\{ 0 \leq \theta \leq \pi/2\}$, one puts a barrier $M(h_2)$ as
in figure \ref{fig:figure10}.

Now do the reflection about the positive $x$
and $y$ axis and then about the negative $x$ and $y$ axes. This gives a complete minimal (Enneper type)
immersed surface. One checks the origin is not a singularity as follows. Reflecting $D_n$ four times
gives an immersed minimal punctured disk that  has the origin in its closure. Gulliver \cite{gulliver}
then proved the origin is not a singularity.

}

\end{document}